\documentclass{article}
\usepackage{amssymb,amsmath,theorem,euscript}

\input{epsf}

\newcounter{sec}

\def\sm{\smallskip}

\newcounter{punct}[sec]

\def\punct{\refstepcounter{punct}{\arabic{sec}.\arabic{punct}.  }}

\def\COUNTERS{\addtocounter{sec}{1}
              \setcounter{punct}{0}
          \setcounter{equation}{0}
          \setcounter{theorem}{0}
                  }
\newtheorem{theorem}{Theorem}[sec]
\newtheorem{proposition}[theorem]{Proposition}
\newtheorem{lemma}[theorem]{Lemma}

\begin{document}

 \def\ov{\overline}
\def\wt{\widetilde}
 \newcommand{\rk}{\mathop {\mathrm {rk}}\nolimits}
\newcommand{\Aut}{\mathop {\mathrm {Aut}}\nolimits}
\newcommand{\Out}{\mathop {\mathrm {Out}}\nolimits}
\renewcommand{\Re}{\mathop {\mathrm {Re}}\nolimits}
\renewcommand{\Im}{\mathop {\mathrm {Im}}\nolimits}

\def\Br{\mathrm {Br}}

\def\SL{\mathrm {SL}}
\def\SU{\mathrm {SU}}
\def\GL{\mathrm {GL}}
\def\U{\mathrm U}
\def\OO{\mathrm O}
 \def\Sp{\mathrm {Sp}}
 \def\SO{\mathrm {SO}}
\def\SOS{\mathrm {SO}^*}
 \def\Diff{\mathrm{Diff}}
 \def\Vect{\mathfrak{Vect}}
\def\PGL{\mathrm {PGL}}
\def\PU{\mathrm {PU}}
\def\PSL{\mathrm {PSL}}
\def\Symp{\mathrm{Symp}}
\def\End{\mathrm{End}}
\def\Mor{\mathrm{Mor}}
\def\Aut{\mathrm{Aut}}
 \def\PB{\mathrm{PB}}
 \def\cA{\mathcal A}
\def\cB{\mathcal B}
\def\cC{\mathcal C}
\def\cD{\mathcal D}
\def\cE{\mathcal E}
\def\cF{\mathcal F}
\def\cG{\mathcal G}
\def\cH{\mathcal H}
\def\cJ{\mathcal J}
\def\cI{\mathcal I}
\def\cK{\mathcal K}
 \def\cL{\mathcal L}
\def\cM{\mathcal M}
\def\cN{\mathcal N}
 \def\cO{\mathcal O}
\def\cP{\mathcal P}
\def\cQ{\mathcal Q}
\def\cR{\mathcal R}
\def\cS{\mathcal S}
\def\cT{\mathcal T}
\def\cU{\mathcal U}
\def\cV{\mathcal V}
 \def\cW{\mathcal W}
\def\cX{\mathcal X}
 \def\cY{\mathcal Y}
 \def\cZ{\mathcal Z}
\def\0{{\ov 0}}
 \def\1{{\ov 1}}
 \def\frA{\mathfrak A}
 \def\frB{\mathfrak B}
\def\frC{\mathfrak C}
\def\frD{\mathfrak D}
\def\frE{\mathfrak E}
\def\frF{\mathfrak F}
\def\frG{\mathfrak G}
\def\frH{\mathfrak H}
\def\frI{\mathfrak I}
 \def\frJ{\mathfrak J}
 \def\frK{\mathfrak K}
 \def\frL{\mathfrak L}
\def\frM{\mathfrak M}
 \def\frN{\mathfrak N} \def\frO{\mathfrak O} \def\frP{\mathfrak P} \def\frQ{\mathfrak Q} \def\frR{\mathfrak R}
 \def\frS{\mathfrak S} \def\frT{\mathfrak T} \def\frU{\mathfrak U} \def\frV{\mathfrak V} \def\frW{\mathfrak W}
 \def\frX{\mathfrak X} \def\frY{\mathfrak Y} \def\frZ{\mathfrak Z} \def\fra{\mathfrak a} \def\frb{\mathfrak b}
 \def\frc{\mathfrak c} \def\frd{\mathfrak d} \def\fre{\mathfrak e} \def\frf{\mathfrak f} \def\frg{\mathfrak g}
 \def\frh{\mathfrak h} \def\fri{\mathfrak i} \def\frj{\mathfrak j} \def\frk{\mathfrak k} \def\frl{\mathfrak l}
 \def\frm{\mathfrak m} \def\frn{\mathfrak n} \def\fro{\mathfrak o} \def\frp{\mathfrak p} \def\frq{\mathfrak q}
 \def\frr{\mathfrak r} \def\frs{\mathfrak s} \def\frt{\mathfrak t} \def\fru{\mathfrak u} \def\frv{\mathfrak v}
 \def\frw{\mathfrak w} \def\frx{\mathfrak x} \def\fry{\mathfrak y} \def\frz{\mathfrak z} \def\frsp{\mathfrak{sp}}
 \def\bfa{\mathbf a} \def\bfb{\mathbf b} \def\bfc{\mathbf c} \def\bfd{\mathbf d} \def\bfe{\mathbf e} \def\bff{\mathbf f}
 \def\bfg{\mathbf g} \def\bfh{\mathbf h} \def\bfi{\mathbf i} \def\bfj{\mathbf j} \def\bfk{\mathbf k} \def\bfl{\mathbf l}
 \def\bfm{\mathbf m} \def\bfn{\mathbf n} \def\bfo{\mathbf o} \def\bfp{\mathbf p} \def\bfq{\mathbf q} \def\bfr{\mathbf r}
 \def\bfs{\mathbf s} \def\bft{\mathbf t} \def\bfu{\mathbf u} \def\bfv{\mathbf v} \def\bfw{\mathbf w} \def\bfx{\mathbf x}
 \def\bfy{\mathbf y} \def\bfz{\mathbf z} \def\bfA{\mathbf A} \def\bfB{\mathbf B} \def\bfC{\mathbf C} \def\bfD{\mathbf D}
 \def\bfE{\mathbf E} \def\bfF{\mathbf F} \def\bfG{\mathbf G} \def\bfH{\mathbf H} \def\bfI{\mathbf I} \def\bfJ{\mathbf J}
 \def\bfK{\mathbf K} \def\bfL{\mathbf L} \def\bfM{\mathbf M} \def\bfN{\mathbf N} \def\bfO{\mathbf O} \def\bfP{\mathbf P}
 \def\bfQ{\mathbf Q} \def\bfR{\mathbf R} \def\bfS{\mathbf S} \def\bfT{\mathbf T} \def\bfU{\mathbf U} \def\bfV{\mathbf V}
 \def\bfW{\mathbf W} \def\bfX{\mathbf X} \def\bfY{\mathbf Y} \def\bfZ{\mathbf Z} \def\bfw{\mathbf w}
 \def\R {{\mathbb R }} \def\C {{\mathbb C }} \def\Z{{\mathbb Z}} \def\H{{\mathbb H}} \def\K{{\mathbb K}}
 \def\N{{\mathbb N}} \def\Q{{\mathbb Q}} \def\A{{\mathbb A}} \def\T{\mathbb T} \def\P{\mathbb P} \def\G{\mathbb G}
 \def\bbA{\mathbb A} \def\bbB{\mathbb B} \def\bbD{\mathbb D} \def\bbE{\mathbb E} \def\bbF{\mathbb F} \def\bbG{\mathbb G}
 \def\bbI{\mathbb I} \def\bbJ{\mathbb J} \def\bbL{\mathbb L} \def\bbM{\mathbb M} \def\bbN{\mathbb N} \def\bbO{\mathbb O}
 \def\bbP{\mathbb P} \def\bbQ{\mathbb Q} \def\bbS{\mathbb S} \def\bbT{\mathbb T} \def\bbU{\mathbb U} \def\bbV{\mathbb V}
 \def\bbW{\mathbb W} \def\bbX{\mathbb X} \def\bbY{\mathbb Y} \def\kappa{\varkappa} \def\epsilon{\varepsilon}
 \def\phi{\varphi} \def\le{\leqslant} \def\ge{\geqslant}

\def\UU{\bbU}
\def\Mat{\mathrm{Mat}}
\def\tto{\rightrightarrows}

\def\Gms{\mathrm {Gms}}
\def\Ams{\mathrm {Ams}}
\def\Isom{\mathrm {Isom}}

\def\Gr{\mathrm{Gr}}

\def\graph{\mathrm{graph}}

\def\O{\mathrm{O}}

\def\la{\langle}
\def\ra{\rangle}


 \def\ov{\overline}
\def\wt{\widetilde}

\renewcommand{\Re}{\mathop {\mathrm {Re}}\nolimits}
\def\Br{\mathrm {Br}}

\def\SL{\mathrm {SL}}
\def\SU{\mathrm {SU}}
\def\GL{\mathrm {GL}}
\def\U{\mathrm U}
\def\OO{\mathrm O}
 \def\Sp{\mathrm {Sp}}
 \def\SO{\mathrm {SO}}
\def\SOS{\mathrm {SO}^*}
 \def\Diff{\mathrm{Diff}}
 \def\Vect{\mathfrak{Vect}}
\def\PGL{\mathrm {PGL}}
\def\PU{\mathrm {PU}}
\def\PSL{\mathrm {PSL}}
\def\Symp{\mathrm{Symp}}
\def\End{\mathrm{End}}
\def\Mor{\mathrm{Mor}}
\def\Aut{\mathrm{Aut}}
 \def\PB{\mathrm{PB}}
 \def\cA{\mathcal A}
\def\cB{\mathcal B}
\def\cC{\mathcal C}
\def\cD{\mathcal D}
\def\cE{\mathcal E}
\def\cF{\mathcal F}
\def\cG{\mathcal G}
\def\cH{\mathcal H}
\def\cJ{\mathcal J}
\def\cI{\mathcal I}
\def\cK{\mathcal K}
 \def\cL{\mathcal L}
\def\cM{\mathcal M}
\def\cN{\mathcal N}
 \def\cO{\mathcal O}
\def\cP{\mathcal P}
\def\cQ{\mathcal Q}
\def\cR{\mathcal R}
\def\cS{\mathcal S}
\def\cT{\mathcal T}
\def\cU{\mathcal U}
\def\cV{\mathcal V}
 \def\cW{\mathcal W}
\def\cX{\mathcal X}
 \def\cY{\mathcal Y}
 \def\cZ{\mathcal Z}
\def\0{{\ov 0}}
 \def\1{{\ov 1}}
 \def\frA{\mathfrak A}
 \def\frB{\mathfrak B}
\def\frC{\mathfrak C}
\def\frD{\mathfrak D}
\def\frE{\mathfrak E}
\def\frF{\mathfrak F}
\def\frG{\mathfrak G}
\def\frH{\mathfrak H}
\def\frI{\mathfrak I}
 \def\frJ{\mathfrak J}
 \def\frK{\mathfrak K}
 \def\frL{\mathfrak L}
\def\frM{\mathfrak M}
 \def\frN{\mathfrak N} \def\frO{\mathfrak O} \def\frP{\mathfrak P} \def\frQ{\mathfrak Q} \def\frR{\mathfrak R}
 \def\frS{\mathfrak S} \def\frT{\mathfrak T} \def\frU{\mathfrak U} \def\frV{\mathfrak V} \def\frW{\mathfrak W}
 \def\frX{\mathfrak X} \def\frY{\mathfrak Y} \def\frZ{\mathfrak Z} \def\fra{\mathfrak a} \def\frb{\mathfrak b}
 \def\frc{\mathfrak c} \def\frd{\mathfrak d} \def\fre{\mathfrak e} \def\frf{\mathfrak f} \def\frg{\mathfrak g}
 \def\frh{\mathfrak h} \def\fri{\mathfrak i} \def\frj{\mathfrak j} \def\frk{\mathfrak k} \def\frl{\mathfrak l}
 \def\frm{\mathfrak m} \def\frn{\mathfrak n} \def\fro{\mathfrak o} \def\frp{\mathfrak p} \def\frq{\mathfrak q}
 \def\frr{\mathfrak r} \def\frs{\mathfrak s} \def\frt{\mathfrak t} \def\fru{\mathfrak u} \def\frv{\mathfrak v}
 \def\frw{\mathfrak w} \def\frx{\mathfrak x} \def\fry{\mathfrak y} \def\frz{\mathfrak z} \def\frsp{\mathfrak{sp}}
 \def\bfa{\mathbf a} \def\bfb{\mathbf b} \def\bfc{\mathbf c} \def\bfd{\mathbf d} \def\bfe{\mathbf e} \def\bff{\mathbf f}
 \def\bfg{\mathbf g} \def\bfh{\mathbf h} \def\bfi{\mathbf i} \def\bfj{\mathbf j} \def\bfk{\mathbf k} \def\bfl{\mathbf l}
 \def\bfm{\mathbf m} \def\bfn{\mathbf n} \def\bfo{\mathbf o} \def\bfp{\mathbf p} \def\bfq{\mathbf q} \def\bfr{\mathbf r}
 \def\bfs{\mathbf s} \def\bft{\mathbf t} \def\bfu{\mathbf u} \def\bfv{\mathbf v} \def\bfw{\mathbf w} \def\bfx{\mathbf x}
 \def\bfy{\mathbf y} \def\bfz{\mathbf z} \def\bfA{\mathbf A} \def\bfB{\mathbf B} \def\bfC{\mathbf C} \def\bfD{\mathbf D}
 \def\bfE{\mathbf E} \def\bfF{\mathbf F} \def\bfG{\mathbf G} \def\bfH{\mathbf H} \def\bfI{\mathbf I} \def\bfJ{\mathbf J}
 \def\bfK{\mathbf K} \def\bfL{\mathbf L} \def\bfM{\mathbf M} \def\bfN{\mathbf N} \def\bfO{\mathbf O} \def\bfP{\mathbf P}
 \def\bfQ{\mathbf Q} \def\bfR{\mathbf R} \def\bfS{\mathbf S} \def\bfT{\mathbf T} \def\bfU{\mathbf U} \def\bfV{\mathbf V}
 \def\bfW{\mathbf W} \def\bfX{\mathbf X} \def\bfY{\mathbf Y} \def\bfZ{\mathbf Z} \def\bfw{\mathbf w}
 \def\R {{\mathbb R }} \def\C {{\mathbb C }} \def\Z{{\mathbb Z}} \def\H{{\mathbb H}} \def\K{{\mathbb K}}
 \def\N{{\mathbb N}} \def\Q{{\mathbb Q}} \def\A{{\mathbb A}} \def\T{\mathbb T} \def\P{\mathbb P} \def\G{\mathbb G}
 \def\bbA{\mathbb A} \def\bbB{\mathbb B} \def\bbD{\mathbb D} \def\bbE{\mathbb E} \def\bbF{\mathbb F} \def\bbG{\mathbb G}
 \def\bbI{\mathbb I} \def\bbJ{\mathbb J} \def\bbL{\mathbb L} \def\bbM{\mathbb M} \def\bbN{\mathbb N} \def\bbO{\mathbb O}
 \def\bbP{\mathbb P} \def\bbQ{\mathbb Q} \def\bbS{\mathbb S} \def\bbT{\mathbb T} \def\bbU{\mathbb U} \def\bbV{\mathbb V}
 \def\bbW{\mathbb W} \def\bbX{\mathbb X} \def\bbY{\mathbb Y} \def\kappa{\varkappa} \def\epsilon{\varepsilon}
 \def\phi{\varphi} \def\le{\leqslant} \def\ge{\geqslant}

\def\UU{\bbU}
\def\Mat{\mathrm{Mat}}
\def\tto{\rightrightarrows}

\def\Gr{\mathrm{Gr}}

\def\graph{\mathrm{graph}}

\def\O{\mathrm{O}}

\def\la{\langle}
\def\ra{\rangle}

\begin{center}
\Large\bf

Biinvariant functions on the group of transformations leaving a measure quasiinvariant
\\
\sc
Neretin Yu.A.%
\footnote{Supported by the grant FWF, P25142.}

\end{center}

{\small Let $\Gms$ be the group of transformations of a Lebesgue space
leaving the measure quasiinvariant, let $\Ams$ be its subgroup consisting of transformations
preserving the measure. We describe canonical forms of double cosets 
of  $\Gms$ by the subgroup  $\Ams$ and show that all
continuous   $\Ams$-biinvariant functions on $\Gms$
are functionals on of the distribution of a Radon--Nikodym derivative.
}

\section{Statements}

\COUNTERS

{\bf\punct The group $\Gms$.%
\label{ss:1}} By $\R^\times$ we denote the multiplicative group of positive reals.
By $t$ we denote the coordinate on $\R^\times$. 

Let $M$ be a Lebesgue space (see \cite{Roh1}) with a continuous probabilistic measure
 $\mu$
(recall that any such space is equivalent to the segment
 $[0,1]$).
Denote by
$\Ams=\Ams(M)$ the group of all transformations (defined up to a.s.)
preserving the measure  $\mu$. By $\Gms=\Gms(M)$ we denote the group of transformations 
(defined up to a.s.) leaving the measure 
  $\mu$ quasiinvariant. 
 
 The group $\Ams$ was widely discussed in connection with ergodic theory,
 the group 
  $\Gms$, which is a topic of the present note, only occasionally was mentioned  in the literature.
However, it is an interesting object from the point of view of representations of infinite-dimensional
groups (``large groups'' in the terminology of A.M.Vershik),
see \cite{Ner-poi}, \cite{Ner-gauss}.

 \smallskip

 {\bf \punct The topology on $\Gms$.}
 A separable topology on $\Gms$ was defined in \cite{Kech} 17.46,
 \cite{Ner-bist}, \cite{Pes},\S 4.5 by different ways. One of the purposes of the present note is two show
 that these ways are equivalent.

 The first way is following.
 Let $A$, $B\subset M$ be measurable subsets.  For $g\in\Gms$
 we define the distribution 
$$\kappa[g;A,B]$$
of the Radon--Nikodym derivative
 $g'$ on the set  
$A\cap g^{-1}(B)$. We say that a sequence 
 $g_j\in\Gms$ converges to  $g$,
if for any measurable sets  $A$, $B$ we have the following weak convergences of measures on $\R^\times$
\begin{equation}
\kappa[g_j;A,B] \to \kappa[g;A,B],\qquad t\kappa[g_j;A,B]\to t\kappa[g;A,B].
\label{eq:convergence}
\end{equation}

{\sc Remark 1.} Point out evident identities:
\begin{equation}
\int_{\R^\times} \kappa[g;A,M](t)=\mu(A),\qquad 
\int_{\R^\times} t\, \kappa[g;A,M](t)=\mu(gA)
\label{eq:id}
.
\end{equation}

{\sc Remark  2.} Consider a measurable finite partition
$$\frh:M=M^1\cup M^2\cup \dots $$ of the space  $M$.
This gives us a matrix  $S_{\alpha\beta}[g;\frh]:=\kappa[g;M^\alpha,M^\beta]$,
composed of measures on   $\R^\times$. If a partition  $\frk$ is a refinement of 
 $\frh$, we  write   $\frh \preccurlyeq\frk$.
Consider a sequence of partitions  $\frh_1 \preccurlyeq \frh_2 \preccurlyeq\dots$,
generating the $\sigma$-algebra of the space%
\footnote{As
$\frh_n$ we can take a partition of the segment 
$M=[0,1]$ into $2^n$ pieces of type  $[k2^{-n}, (k+1)2^{-n})$} $M$.
 A convergence $g_j\to g$ is equivalent to an element-wise convergence
 in the sense 
  (\ref{eq:convergence})
of all matrices $S[g_j;\frh_n]\to S[g;\frh_n]$.

\begin{proposition}
\label{pr:1}
The group $\Gms$ is a Polish group with respect to this topology,
i.e., $\Gms$ is a  separable topological group  complete with respect to the two-side uniform structure and
homeomorphic to a complete metric space%
\footnote{A metric is compatible with the topology of
the group, but not with its algebraic structure; in particular a metric is not assumed to be invariant.
A completeness of a group  in the sense of two-side uniform structure
 (in Raikov's sense \cite{Rai}) is defined (for metrizable groups)
 in the following way. Let double sequences  
 $g_ig_j^{-1}$ and $g_i^{-1}g_j$ converge to  1 as $i$, $j\to\infty$. Then $g_i$
 has a limit in the group.This definition is not equivalent to the definition
 of Bourbaki  \cite{Bou}, III.3.3, who requires a completeness with respect to
 both one-side uniform structures. The group  $\Gms$ is not complete in the sense of Bourbaki.}. 
\end{proposition}

Let $1\le p\le \infty$, $s\in\R$. The group $\Gms$ acts in the space  $L^p(M)$
by isometric transformations  according the formula 
$$
T_{1/p+is}f(x)=f(g(x))g'(x)^{1/p+is}
.
$$
On the space
 $\cB(V)$ of operators of a Banach space   $V$ we define in the usual way
 (see, e.g., \cite{RS}, VI.1) the strong and weak topologies.  
 Also, on the set $\cG\cL(V)$ of invertible operators
 we introduce a {\it bi-strong} topology, $A_j$ converges to $A$, 
if $A_j\to A$ and $A^{-1}_j\to A^{-1}$ strongly. The embedding $T_{1/p+is}:\Gms \to \cB(L^p)$
induces a certain topology on
 $\Gms$ from any operator topology on  $\cB(V)$ or $\cG\cL(V)$.

\begin{proposition}
\label{pr:2}
{\rm a)} Let $1<p<\infty$, $s\in\R$. A topology on  $\Gms$ induced from any of three
topologies  {\rm(}strong, weak, bi-strong{\rm)} coincides with the topology
defined above.

\sm

{\rm b)} Let $p=1$, $s\in\R$. A topology on $\Gms$ induced from strong or bi-strong
topology coincides with the topology defined above.

\sm

{\rm c)} Let  $1\le p<\infty$, $s\in\R$. Then the image of $\Gms$ в $\cG\cL(L^p(M))$ is closed 
in the bi-strong topology.
\end{proposition} 

Point out that the coincidence of topologies is not surprising. 
It is known that two different Polish topologies on a group can not determine 
the same Borel structure,
 see \cite{Kech}, 12.24. There are also theorems
 about automatic continuity of homomorphisms, see  \cite{Kech},  9.10,

\sm

{\bf\punct Double cosets $\Ams\setminus \Gms/\Ams$.
Canonical forms.%
\label{ss:predst}} We reformulate the problem of description of double cosets
 $\Ams\setminus \Gms/\Ams$ in the following way. Let  
$(P,\pi)$, $(R,\rho)$ be Lebesgue spaces with continuous probabilistic
measures. Denote by 
 $\Gms(P,R)$ the space of all bijections  $g:P\to R$
 (defined up to a.s.), such that images and preimages of sets of zero measure have zero measure.
We wish to describe such bijections up to the equivalence 
  \begin{equation}
  g\sim u\cdot g\cdot v,\qquad \text{where $v\in\Ams(P)$, $u\in \Ams(R)$}
  \label{eq:sim}
  \end{equation}
  (clearly, such classes are in-to-one correspondence with double cosets
 $\Ams\setminus\Gms/\Ams$).
  
  \begin{lemma}
  \label{l:derivative}
Two elements  
  $g_1$, $g_2\in\Gms(P,R)$ are contained in one class if and only if the
Radon--Nikodym derivatives
 $g_1'$, $g'_2:P\to \R$
are equivalent with respect to the action of the group
 $\Ams(P)$, i.e.,
  $g_2'(m)=g_1'(hm)$, where $h$ is an element of  $\Ams(P)$.
  \end{lemma}
  
An evident invariant of this action is the distribution 
 $\nu$ of the Radon--Nikodym derivative
$g'$
of the map $g$,
\begin{equation}
\int_{\R^\times} d\nu(t)=1 ,\qquad \int_{\R^\times} t\,d\nu(t)=1.
\label{eq:nu}
\end{equation}
This invariant is not  	exhaust, the problem is reduced to the Rokhlin theorem 
 \cite{Roh2} on metric classification of functions, see discussion below, \S
 \ref{ss:roh}. The final answer is following. 

 Consider a countable number of copies 
 $\R_1^\times$, $\R_2^\times$, \dots of half-line  $\R^\times$.
 Consider one more more copy  $\R_\infty^\times$. Consider the disjoint union
$$
\cL:=\R_1^\times\coprod \R_2^\times\coprod \R_3^\times\coprod \dots \coprod \bigl(\R_\infty^\times \times [0,1]\bigr).
$$
Let $\nu_1$, $\nu_2$,\dots и $\nu_\infty$ be a family of measure on  $\R^\times$
satisfying the following conditions

1. $\nu_1$, $\nu_2$, \dots are continuous (but $\nu_\infty$ admits atoms).

2. $\nu_1\ge\nu_2\ge\dots$

3. The measure   $\nu:=\nu_1+\nu_2+\dots+\nu_\infty$ satisfies (\ref{eq:nu}).

Equip each  $\R^\times_j$ with the measure  $\nu_j$, equip $\R^\times_\infty\times[0,1]$ with
the measure $\nu_\infty\times dx$, where $dx$ is the Lebesgue measure on the segment. 
Denote the resulting measure space by $\cL[\nu_1,\nu_2,\dots;\nu_\infty]$.

Consider the same measure on  $\cL$ multiplied by   $t$, we denote the resulting measure space
by $\cL_*[\nu_1,\nu_2,\dots;\nu_\infty]$.

Consider the identity map
\begin{equation}
{\rm id}:\cL[\nu_1,\nu_2,\dots;\nu_\infty]\to \cL_*[\nu_1,\nu_2,\dots;\nu_\infty]
\label{eq:identity}
\end{equation}
Evidently, the distribution of the Radon--Nikodym derivative 
of the map ${\rm id}$ coincides with $\nu$.

\begin{proposition}
\label{pr:3}
Any equivalence class 
 {\rm(\ref{eq:sim})} contains a unique representative of the type 
  {\rm(\ref{eq:identity})}.
\end{proposition} 

Denote the double coset containing this representative by
$S[\nu_1,\nu_2,\dots;\nu_\infty]$.

\sm

{\bf\punct On closures of double cosets.}


\begin{theorem}
\label{l:4} Let a measure  $\nu$ on $\R^\times$ satisfies {\rm(\ref{eq:nu})}, let
$\nu=\nu^c+\nu^d$ be its decomposition into continuous and discrete parts.
Then the closure of the double coset
$S[\nu^c,0,0,\dots;\nu^d]$
contains all double cosets  
$S[\nu_1,\nu_2,\nu_3,\dots;\nu^c_\infty+\nu_d]$ with $\nu_1+\nu_2+\dots+\nu_\infty^c=\nu^c$.
\end{theorem}

{\bf\punct Hausdorff quotient.}
Consider the space $\cM$ of all measures  $\nu$
on $\R^\times$ satisfying  (\ref{eq:nu}). Say that  $\nu^j\in\cM$ converges
to $\nu$  if $\nu^j\to \nu$ and $t\nu^j\to t\nu$ weakly.

Consider a map
$\Phi:\Gms\to \cM$  that for any  $g$ assigns the distribution
of its Radon--Nikodym derivative
(i.e., $\Phi(g)=\kappa[g;M,M]$). In virtue of Theorem   \ref{l:4},
preimages of points  $\nu\in\cM$ are closures of double cosets
$S[\nu^c,0,0,\dots;\nu^d]$.



\begin{theorem}
\label{th:quotient}
Let  $f$ be a continuous map of $\Gms$ to a metric space  $T$,
moreover, $f$ let be constant on double cosets.
Then  $f$ has the form  $f=q\circ \Phi$, where $q:\cM\to T$ is a continuous map.
\end{theorem}

\sm

{\bf\punct A continuous section $\cM\to\Gms$.}
We say that a function  $h:[0,1]\to[0,1]$ is contained in the class  $\cG$,
if

$\bullet$ $h$ is downward convex;

$\bullet$ $h(0)=0$, $h(1)=1$, and $h(x)>0$ for $x>0$.

Any such function is an element of the group
$\Gms\bigl([0,1]\bigr)$.

\begin{proposition}
\label{pr:section}
Let $\nu\in \cM$. Then there is a unique function 
$\psi:[0,1]\to[0,1]$ of the class  $\cG$ such that the distribution
of the derivative $\psi'$
is  $\nu$. Moreover, the map  $\nu\mapsto\psi$ is a continuous map  $\cM\to\Gms$.
\end{proposition}

{\bf\punct A more general statement.}
Consider a finite or countable measurable partition of our measure space 
$M=\coprod_j M_j.$
Denote by 
 $K$ the direct product 
$
K=\Ams(M_1)\times \Ams(M_2)\times \dots 
$.
Consider the double cosets 
 $K\setminus \Gms/K$.
Assign to each  $g\in \Gms$ the matrix $\kappa_{ij}=\kappa[g;M_i,M_j]$
composed of measures on  $\R^\times$. Denote by  $\cS$ the set of matrices
that can be obtained in this way, i.e.,
$$
\sum_{j} \int_{\R^\times} d\kappa_{ij}(t)=\mu(M_i),
\qquad
\sum_{i} \int_{\R^\times} t\kappa \, d\mu_{ij}(t)=\mu(M_j).
$$
Equip  $\cS$ with element-wise convergence  (\ref{eq:convergence}).
Denote by  $\Psi$ the natural map  $\Gms\to \cS$.

\begin{theorem}
\label{th:last}
 Let  $f$ be a continuous map from $\Gms$ to a metric space 
 $T$. Then there exists a continuous map  $q:\cS\to T$,
 such that $f=q\circ \Psi$.
\end{theorem}

Point out that this statement was actually used in 
 \cite{Ner-bist}, \cite{Ner-poi}.

\sm

{\bf \punct The structure of the note.}
The statements about topology on $\Gms$ are proved in
 \S2, about double cosets in   \S3.  Theorem \ref{th:quotient} follows from
 Theorem
\ref{l:4}. However, as the referee pointed out, the first statement
is simpler than the second
(and it is more important). Therefore in the beginning of 
  \S3 we present a separate proof of Theorem  \ref{th:quotient}.




\section{The topology on the group $\Gms$}

\COUNTERS

Below we prove Propositions \ref{pr:1} and \ref{pr:2}. The main 
	auxiliary statement is Lemma 
\ref{l:2.4}. The remaining lemmas are proved in a straightforward way.

\smallskip

Notation:

$\bullet$ $\delta_a$ is an probabilistic atomic measure
 $\R^\times$ supported by a point  $a$.

$\bullet$  $\{\cdot,\cdot\}_{pq}$ is the natural pairing of $L^p$ and $L^q$, where
$1/p+1/q=1$;

$\bullet$  $\chi_A$ is the indicator function of a set  $A\subset M$,
i.e., $\chi_A(x)=1$ for $x\in A$ and $\chi_A(x)=0$ for $x\notin A$.

\sm

{\bf\punct  Preliminary remarks on the spaces  $L^p$.}

\sm

1)
Recall  (see \cite{FJ}, \S
3.3) that for  $p\ne 2$ the group of isometries  $\Isom\bigl(L^p(M)\bigr)$
of the space   $L^p(M)$ consists of operators of the form 
\begin{equation}
R(g, \sigma)f(x)=\sigma(x) f(g(x))g'(x)^{1/p}
,
\label{eq:isometries}
\end{equation}
where $g\in\Gms$, and  $\sigma:M\to\C$ is a function whose absolute value equals 1.

2)
For  $1<p<\infty$ the space   $L^p$ is uniformly convex (see \cite{Koth}, \S 26.7), 
therefore the restrictions of the strong and weak topologies to the unit sphere coincide.
Therefore on the group of isometries  $\Isom\bigl(L^p(M)\bigr)$
the weak and strong operator topologies coincide.

3) Recall that for separable Banach spaces 
(in particular, for  $L^p$ with $p\ne \infty$) the group of all isometries
equipped with bi-strong topology is a Polish group, see  \cite{Kech}, 9.B9.

\sm

{\bf\punct Preliminary remarks on the group  $\Gms$.}

\sm

1) {\it The invariance of the topology.} Equip   $\Gms$ with topology from Subsection 1.2.
The product in $\Gms$ is separately continuous (this is a special case of Theorem
 5.9 
from \cite{Ner-boundary}). In particular, this implies that the topology
on $\Gms$ is invariant with respect to left and right shifts.

 The map $g\mapsto g^{-1}$ is continuous. Indeed,
 $$
 \kappa[g^{-1};B,A](t)=t^{-1}\kappa[g;A,B](t^{-1})
 ,$$
 and this map transpose the convergences 
 (\ref{eq:convergence}).

\sm

2) {\it Separability of $\Gms$}. For a measure  $\kappa[g;A,B]$ 
consider the {\it characteristic function}
\begin{equation}
\chi(z)=\int_{\R^\times} t^z d \kappa[g;A,B](t)
,
\label{eq:char}
\end{equation}
continuous in the strip  $0\le  \Re z\le 1$ and holomorphic in the open strip.
The convergence of measures  $\kappa$ is equivalent to point-wise
convergence of characteristic functions
uniform in each rectangle
$$0\le  \Re z\le 1,\qquad -N\le \Im z\le N, $$
\cite{Ner-boundary}, Propositions 4.4-4.5. This convergence is 
separable. Next, by Remark 2 of \S\ref{ss:1},
it suffices to verify 
the convergence of measures  $\kappa[g_j;A,B]\to \kappa[g;A,B]$ 
for an appropriate countable set of pairs measurable subsets
$(A,B)$.

\sm

3) {\it The action on Boolean algebra of sets.}

\begin{lemma}
\label{eq:bigtriangleup}
Let $g_j\to g$ in $\Gms$. Then for any measurable set $A\subset M$
we have
\begin{equation}
 \mu(g_j A\bigtriangleup gA)\to 0
 \label{eq:triangle}
.\end{equation}
\end{lemma}

{\sc Proof.} By the invariance of the topology it suffices 
to consider $g=1$. Then
$$
\mu(g_j A\cap A)=\int\limits_{\R^\times} d\kappa[g_j^{-1}; A;  A](t)
\to
\int\limits_{\R^\times} d\kappa[1; A;  A](t)=
\int\limits_{\R^\times} \mu(A)\delta_0(t)=\mu(A)
;$$
$$
\mu(g_j A)=\int\limits_{\R^\times} t\,d\kappa[g_j;A;M]\to 
\int\limits_{\R^\times} t\,d\kappa[1;A;M]=\int\limits_{\R^\times} \mu(A) t\,\delta_0(t)=\mu(A)
.
$$
Comparing two rows we get the desired statement.
\hfill $\square$

\sm

{\sc Remark.} The opposite is false. Let $M=[0,1]$, 
$$
g_j(x)=x+\frac 1{2\pi n}\sin(2\pi n x).
$$
Then for any
 $A\subset[0,1]$ we have   $\mu(g_j(A)\bigtriangleup A)\to \mu(A)$.
But there is no convergence $g_j\to 1$ in $\Gms$; $T_1(g_j)$ converges weakly
to 1 in $L^1$, but there is no strong convergence. 
\hfill $\boxtimes$

\sm

4){\it  The continuity of representations $T_{1/p+is}$.}

\begin{lemma}
For $p<\infty$ the homomorphisms  $T_{1/p+is}:\Gms\to \Isom(L^p)$
are continuous with respect to the weak topology
 $\Isom(L^p)$.
\end{lemma}

{\sc Proof.}
Let $g_j\to g$ in $\Gms$. Consider 'matrix elements'  
$$
\{ T_{1/p+is}(g_j)\chi_A,\chi_B\}_{pq}=
\int\limits_{A\cap g_j^{-1}B} g_j'(x)^{1/p+is}d\mu(x)=
\int\limits_{\R^\times}t^{1/p+is}\,d\kappa[g_j;A,B](t)
$$
Weak convergence of measures
(\ref{eq:convergence}) implies the convergence of characteristic functions 
(\ref{eq:char}), our expression tends to  
$$
\int\limits_{\R^\times}t^{1/p+is}\,d\kappa[g;A,B](t)=
\int\limits_{A\cap g^{-1}B} g'(x)^{1/p+is}d\mu(x)=\{ T_{1/p+is}(g)\chi_A,\chi_B\}_{pq}
,$$
as required.
\hfill $\square$

\sm

Thus, for $1<p<\infty$ the maps  $T_{1/p+is}:\Gms\to\Isom(L^p)$  are continuous with respect
to the strong (=weak) topology. Keeping in mind the continuity
of the map $g\mapsto g^{-1}$, we get that the maps $T_{1/p+is}$ are continuous with respect to the bi-strong
topology.

The case $L^1$ must be considered separately. 

\begin{lemma}
Let $g_j\to g$. Then $T_{1+is} (g_j)\in \Isom(L^1)$ strongly converges to 
$T_{1+is} (g)$.
\end{lemma}

{\sc Proof.} Without loss of generality, we can set   $g=1$.
It suffices to verify the convergence   $\|T_{1+is} (g_j)\chi_A-\chi_A\|\to0$
for any measurable  $A$. This equals
\begin{multline}
\int_M \bigl|\chi_A(g_j x) g'(x)^{1+is}-\chi_A(x)\bigr|\,d\mu(x)=\\=
\int_{A\cap g_j^{-1}A}\bigl|  g'(x)^{1+is}-1\bigr|\,d\mu(x)+
\int_{A\setminus g_j^{-1}A}d\mu(x)+ \int_{ g_j^{-1}A\setminus A}  g'(x)\,d\mu(x)
=\\=
\int_{\R^\times}|t^{1+is}-1|\,d\kappa[g_j;A,A](t)+
\mu\bigl(A\setminus g_j^{-1}A\bigr)+\mu(A\setminus g_j A)
\label{eq:L1}
.
\end{multline}
The second and the third summands tend to 0 by Lemma
 \ref{eq:bigtriangleup},
measures  $\kappa[\dots]$ and $t\kappa[\dots]$ converge weakly to
$\mu(A)\delta_0$, therefore the integral tends to 0. 
\hfill $\square$

\sm

{\bf\punct The coincidence of topologies and the continuity of the multiplication.}

\begin{lemma}
\label{l:2.4}
Let $1<p<\infty$.
Let $T_{1/p+is}(g_j)$ weakly converge to   $1$ in $\Isom(L_p)$. 
Then $g_j$ converges to  $1$ in $\Gms$.
\end{lemma}

{\sc Proof.}
{\it Step 1.} Now it will be proved that  $g'_j$ converges to   1 in $L^1(M)$.
For this purpose, we notice that the following
sequence of matrix elements must converge to 
 1:
\begin{equation}
\{T_{1/p+is}(g_j)\, 1, 1\}_{pq}=\int\limits_M g'_j(x)^{1/p+is}d\mu(x)
=\int\limits_{\R^\times} t^{1/p+is} d\kappa[g_j;M,M](t).
\label{eq:me}
\end{equation}
Estimate the integrand:
$$
\Re t^{1/p+is}\le  t^{1/p}\le \frac 1 q+ \frac t p.
$$
The second inequality means that the graph of upward convex function is lower than
the tangent line at 
 $t=1$. From another hand:
\begin{multline*}
\int_{\R^\times}\Bigl( \frac 1 q+ \frac t p\Bigr)\, d\kappa[g_j;M,M](t)
=\\
=\frac 1 q\int_{\R^\times}d\kappa[g_j;M,M](t)+ 
\frac 1 p\int_{\R^\times}t\,d\kappa[g_j;M,M](t)=\frac 1 q+ \frac 1 p=1.
\end{multline*}
Look to a deviation of integral 
 (\ref{eq:me}) from 1.
 The same reasoning with tangent line allows to estimate the difference
$
\frac 1 q+ \frac t p-t^{1/p}
$.
For any
  $\epsilon>0$ there is   $\sigma>0$ such that
$$
\frac 1 q+ \frac t p-t^{1/p}>\begin{cases}
\sigma\qquad \text{for $t<1-\epsilon$};\\
\sigma t\qquad\text{for $t>1+\epsilon$}.
\end{cases}
$$
Therefore
$$
1-\Re \{T_{1/p+is}(g_j)\, 1, 1\}_{pq}>\sigma\int\limits_0^{1-\epsilon}d\kappa[g_j;M,M](t)+
\sigma
\int\limits_{1+\epsilon}^\infty t
\, d\kappa[g_j;M,M](t)
.
$$
This must tend to 0, therefore
$\kappa[g_j;M,M]$ and $t\cdot \kappa[g_j;M,M]$
tend to $\delta_0$ weakly. This implies the convergence 
$g'_j\to 1$ in the sense of $L^1$.

The remaining part of the proof is more-or-less automatic.

\sm

{\it Step 2.} 
Let  $z$ be contained in the strip $0\le \Re z\le 1$.
Let us show that  $(g'_j)^{z}$  tends to 1  in the sense of  $L^1$.
Let  $\|g'-1\|_{L^1(M)}<\epsilon$. Then there is an uniform with respect to
 $g$
estimate $\|(g')^z-1\|_{L^1(M)}<\psi_z(\epsilon)$, where  $\psi_z(\epsilon)$ tends to  0
as  $\epsilon$ tends to 0. For this aim it is sufficient to notice that
$$
|a^z-1|< \begin{cases}
|z|\,(a-1) \qquad &\text{for $a>1$};\\
|z|\, 2^{-\Re z+1} |a-1| \qquad  &\text{for $1/2\le a\le 1$};\\
2 \qquad \text{for $0<a<1/2$},
\end{cases}
$$
moreover,  $g'<1/2$ can be  only on the set of measure $\le 2\epsilon$.

\sm

In particular, for any subset
$C\subset M$ we have 
\begin{equation}
\Bigl|\int_C g'(x)^z dx-\mu(C)\Bigr|\le\psi_z(\epsilon)
\label{eq:psi()}
.
\end{equation}

\sm

{\sc Step 3.} Now we use convergence of matrix elements:
$$\{ T_{1/p+is}(g_j)\chi_A,\chi_B\}_{pq}
=\int\limits_{A\cap g_j^{-1}B} g_j'(x)^{1/p+is}d\mu(x)
\to 
\{ \chi_A,\chi_B\}_{pq}= \mu(A\cap B).$$
By (\ref{eq:psi()}),  we have convergence
$$
\int_{A\cap g_j^{-1}B} g_j'(x)^{1/p+is}d\mu(x)-\mu(A\cap g_j^{-1}B)\to 0
.
$$
Comparing two last convergences we get
$\mu(A\cap g_j^{-1}B)
\to \mu(A\cap B)$.

\sm

{\it Step 4.} By the convergence  $(g_j')^z$ in $L^1(M)$, we have
$$
\int t^z \,d\kappa[g_j;A,B](t)=\int_{A\cap g_j^{-1}B}g'(x)^z\,dx \to
\mu(A\cap B)
$$
for each $z$; the point-wise convergence of characteristic functions
implies weak converges (\ref{eq:convergence}) of measures
(see, \cite{Ner-boundary}), in our case, to $\mu(A\cap B)\delta_0$.
\hfill $\square$

\sm

Thus the topology on
$\Gms$ is induced from the strong operator topology of the spaces
 $L^p$. In separable Banach spaces the multiplication
is continuous in the strong topology on bounded sets.
Therefore, the multiplication in  $\Gms$ is continuous.

\begin{lemma}
Let operators $T_{1+is}(g_j)$ converge to $1$ in the strong operator topology
of spaces $L^1$.
Then  $g_j\to 1$ in $\Gms$.
\end{lemma}

{\sc Proof.} In  (\ref{eq:L1}) the first row must tend to zero.
Therefore all summands of the last row tend 0, in particular the first one.
This implies weak convergences of measures   $\kappa[g_j;A,A]$ and $t\kappa[g_j;A,A]$ to
$\mu(A)\delta_1$. Comparing this with  (\ref{eq:id}), we get  convergences
$\kappa[g_j;A,M\setminus A]$ and $t\kappa[g_j;A,M\setminus A]$
to 0. Now it is easy to derive the convergence of  $g_j\to g$ in $\Gms$.
\hfill $\square$

\sm

{\bf \punct The completeness of $\Gms$.}
The group of isometries of a separable Banach space is a Polish group with respect to
the bi-strong topology
 (\cite{Kech}, 9.3.9). Let $p\ne 1$
$2$, $\infty$ и $s=0$. Then the isometries  $T_{1/p}(g)$ are precisely isometries
(\ref{eq:isometries}) that send the cone of non-negative functions 
to itself. 
Obviously, the set of operators sending this cone to itself is weakly closed.
Therefore, 
 $\Gms$ is a closed subgroup in the group of all isometries
and therefore it is complete. 

\sm

\begin{figure}
$$
\epsfbox{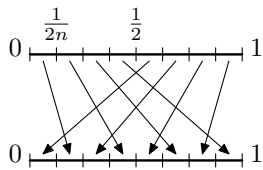}
$$
\caption{Reference to Example \ref{ss:closure}.\label{fig:1}}
\end{figure}

{\bf \punct Bi-strong closeness of the image.%
\label{ss:closure}} The group $\Gms$ is closed in the group  $\Isom(L^p)$,
since it is complete with respect of the induced topology.

\sm

It is noteworthy that the group 
$\Isom(L^p)$ is not strongly closed in the space of bounded operators
in  $L^p$. 
The images of the  groups $\Ams$ and $\Gms$ also are not closed.

\sm

{\sc Example.} Let $p\ne\infty$. Consider an operator in $L^p$ of the form 
$$
Rf(x)=
\begin{cases}
f(2x), \qquad\text{ for $0\le x\le 1/2$};
\\
f(2x-1), \qquad\text{for $1/2< x\le 1$};
\end{cases}
$$
For any function
 $f$ we have   
$\|Rf\|=\|f\|$. However, this operator is not invertible.
For the sequence $g_n\in \Ams$  from Fig.
 \ref{fig:1}
we have the strong convergence  $T_{1/p}(g_n)$ to $R$.
\hfill $\boxtimes$.

\sm

Weak closures for some subgroups
 $\Gms$ are discussed in 
 \cite{Ner-poi}, \cite{Ner-gauss}.

\section{Double cosets}

\COUNTERS

{\bf\punct Proof of Theorem \ref{th:quotient}.}
Denote by
 $G^0\subset\Gms$ the group of transformations whose Radon--Nikodym derivative
 has only finite number of values. Obviously,

\sm

$\bullet$ The subgroup $G^0$ is dense in $\Gms$.

\sm

$\bullet$  Double cosets  $\Ams\setminus G^0/\Ams$ are completely determined
by the distribution of the Radon--Nikodym derivative.

\sm

Consider a measure  $\kappa\in\cM$. Consider a sequence of discrete measures 
$\kappa_N\in \cM$ convergent to $\kappa$ and having the following property:
Fix  $N$ and cut the semi-axis $t>0$ into pieces of length  $2^{-N}$.
For any  $j\in \N$ we require the following coincidence of measures
of semi-intervals
$$
\int\limits_{
\frac {j-1}{2^N}<t\le  \frac{j}{2^N}}\!\!
d\kappa(t)=
\!\!
\int\limits_{
\frac {j-1}{2^N}<t\le  \frac{j}{2^N}}\!\!
d\kappa_N(t),\qquad
\int\limits_{
\frac {j-1}{2^N}<t\le  \frac{j}{2^N}}
\!\!
t\cdot
d\kappa(t)=
\!\!
\int\limits_{
\frac {j-1}{2^N}<t\le  \frac{j}{2^N}}
\!\!
t\cdot
d\kappa_N(t)
$$

Consider  $g\in \Gms$ whose distribution of the Radon--Nikodym derivative equals $\kappa$.
Consider a sequence  $g_N\in G^0$ convergent to $g$ such that 
a distribution of the Radon--Nikodym derivative of 
$g_N$ is $\kappa_N$. For this, we fix  
$N$ and for each  $j$
consider the subset  $A_j\subset M$, where the Radon--Nikodym derivative
satisfies 
$$\frac {j-1}{2^N}<g'(x)\le \frac{j}{2^N}.
$$
Set $B_j=g(A)$. 
Consider an arbitrary map 
 $g_N\in G^0$ such that  
$g_N$ send  $A_j$ to $ B_j$ and
the distribution of the Radon--Nikodym derivative of
 $g_N$ coincides with the restriction of  
the measure $\kappa_N$ of the semi-interval 
$\bigl(\frac {j-1}{2^N}<t\le \frac{j}{2^N}]$. It easy to see that the sequence
$g_N$ converges to  $g$.

Now, let $f$ be a continuous function on $\Gms$ constant on double cosets.
Пусть $g$ and  $h\in\Gms$ have same distribution of Radon--Nikodym derivatives.
Then 
$g_N$ and $h_N$ are contained in the same double coset, wherefore $f(g_N)=f(h_N)$.
By continuity of  $f$ we get $f(g)=f(h)$. 

To avoid a proof of the continuity the map
 $q$ (see the statement of the theorem), we refer to Proposition 
\ref{pr:section} (which is proved below independently 
of the previous considerations).

\sm

{\bf\punct Proof of Proposition \ref{pr:3}.%
\label{ss:roh}}
Let  $M\simeq[0,1]$ be a Lebesgue space.
Invariants of measurable functions  $f:M\to \R$
with respect to the action of $\Ams(M)$ were described by Rokhlin
in  \cite{Roh2}.
To any function $f$ he assigns its distribution function  $F(y)$, i.e., the measure 
of the set
$M_y\subset M$ determined by the inequality $f(x)< y$. Also he assigns to $f$ 
a sequence of functions
$F_1$, $F_2$, \dots, where $F_n(y)$ is the supremum
of measures of all sets  $A\subset M_y$, 
on which $f$ takes each value   $\le n$ times. These data satisfy the following conditions:

$\bullet$ the function $F$ satisfies the usual properties of 
distribution functions:
$F$ is a left-continuous non-decreasing function, $\lim_{y\to-\infty} f(y)=0$,
$\lim_{y\to+\infty} f(y)=1$;

$\bullet$ $F_n$  are non-decreasing functions;

$\bullet$
$0\le F_1(y)\le F_2(y)\le\dots\le F(y)$;

 $\bullet$ $F_k(y)-2 F_{k+1}(y)+ F_{k+2}(y)\ge 0$ for all $k$.

 \sm

According  \cite{Roh2}, a function $f$  determined up to the action of the group
 $\Ams$ is uniquely defined by the invariants
 $F_1$, $F_2$,\dots, $F$. Moreover, for any collection of functions
 $F_1$, $F_2$,\dots, $F$ with above listed properties there exists
 $f$, whose invariants coincide with  $F_1$, $F_2$,\dots, $F$.

Now we will describe canonical forms of functions
 $f$ under the action of the group  $\Ams$.
Consider a collection of continuous measures 
$\nu_1\le\nu_2\le\dots$ on $\R$ and the measure $\nu_\infty$ on $\R$
such that 
$\nu_1(\R)+\nu_2(\R)+\dots+\nu_\infty(\R)=1$. Denote by   $t$ the coordinate on 
$\R$. Consider the disjoint union of the spaces with measures
\begin{equation}
\cL=\Bigl((\R,\nu_1)\coprod(\R,\nu_2)\coprod\dots\Bigr)\coprod (\R\times[0,1],\nu_\infty \times ds)
,
\label{eq:for-roh}
\end{equation}
where $ds$ is the Lebesgue measure on the segment  $[0,1]$.
Consider the function  $f$ on $\cL$ that equals to  $t$ 
on each copy of $\R$ and  equals to  $t$ on $\R\times[0,1]$. 

The invariants of this function are
$$
F_n(y)=\sum_{j\le n}\nu_j(-\infty,y),\qquad F(y)=\sum_{1\le j<\infty} \nu_j(-\infty,y)+
\nu_\infty(-\infty,y)
$$
It can be readily seen that measures
 $\nu_1$, $\nu_2$,\dots, $\nu_\infty$ admit a reconstruction from the invariants
$F_1$, $F_2$,\dots, $F$. Moreover any admissible collection of invariants
corresponds to a certain collection of measures  
$\nu_1$, $\nu_2$, \dots, $\nu_\infty$.

Now consider an element 
 $g\in \Gms(P,R)$. Reduce the derivative  $g':P\to\R^\times$   
 to the canonical form by a multiplication
 $g\mapsto gh$, where $h\in\Ams$. 
Since $g'(x)>0$, all the measures  
$\nu_j$, $\nu$ are supported by the half-line $t>0$. The integral of $g'$   is 1,
therefore
\begin{equation}
\sum_j\int t\,d\nu_j(t)+\int t\,d\nu_\infty(t)=1.
\label{eq:for-prob}
\end{equation}
Now we assume
 $P=\cL$, see (\ref{eq:for-roh}). Let $\cL_*$ be obtained from   $\cL$ by a multiplication 
 of the measure by $t$.
In virtue of  (\ref{eq:for-prob}), this measure must be probabilistic.
The map $g:\cL\to R$ can be regarded as a map 
 $g_*:\cL_*\to R$. Since  $g'=t$, for any measurable set  $B\subset\cL$ the measure of $B$
in $\cL_*$ coincides with the measure  $g(B)$. Therefore $g_*:\cL_*\to R$ preserves measure.

Thus $g$ is reduced to the canonical form.

 
 \sm

{\bf\punct Splitting of measures.} We start a proof of Theorem   \ref{l:4}.
Modify the notation for  $\cL[\nu_1,\nu_2,\dots;\nu_\infty]$, 
$\cL_*[\nu_1,\nu_2,\dots;\nu_\infty]$ and   $S[\nu_1,\nu_2,\dots;\nu_\infty]$ from (\ref{ss:predst}).
Now it is convenient to reject the condition  
 $\nu_1\ge\nu_2\ge\dots$. Also, we weaken condition 
(\ref{eq:nu}) and set
\begin{equation}
\int_{\R^\times} d\nu(t)<\infty ,\qquad \int_{\R^\times} t\,d\nu(t)<\infty.
\label{eq:nu-1}
\end{equation}

 Let $\nu$ be a continuous measure  on $\R^\times$ satisfying 
(\ref{eq:nu-1}). Consider the space  
$\cL[\nu, 0,0,\dots;0]$. Represent $\nu$ as a sum $\nu=\nu_1+\nu_2$.

\begin{lemma}
\label{l:split}
The closure of the class  $S[\nu, 0,0,\dots;0]$ contains 
$S[\nu_1, \nu_2,0,\dots;0]$.
\end{lemma}

\begin{figure}
$$
\epsfbox{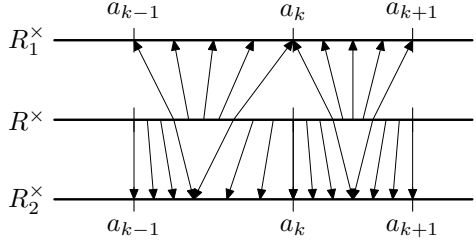}
$$
\caption{A reference to Lemma  3.1. The map[ $\phi_n$. \label{fig:2}}

\end{figure}

{\sc Proof.}  Denote
$$
\cL:=\cL[\nu_1,0,\dots;0],\qquad \cL':=\cL[\nu_1, \nu_2,0,\dots;0]
.
$$
The same measure spaces with the measure multiplied by 
 $t$ we denote as
$\cL_*$,  $\cL_*'$
 Now we will construct two sequences of measure preserving bijections
$$\phi_n: \cL\to \cL', 
\qquad
\psi_n: \cL_*\to \cL'_*
.$$
Cut  $(\R^\times,\nu)$ by $2^n$ intervals $C_0$,\dots,$C_{2^n-1}$ by points
$$
a_k=\frac{k 2^{-n}}{1-k 2^{-n}} ,\qquad k=1,2,\dots,2^n-1
.$$
Denote this partition%
\footnote{The only necessary for us  property of partition is the following: a diameter of a partition
on any finite interval  $(0,M]$ tends to 0 as 
 $n\to\infty$.}  by $\frh_n$.

\begin{lemma}
The exists a sub-interval 
 $B_k\subset C_k$ such that $\nu(B_k)=\nu_1(C_k)$,
$(t\cdot\nu)(B_k)=(t\cdot\nu_1)(C_k)$.
 \end{lemma}

{\sc Proof.} We have $C_k=[a_k,a_{k+1}]$. Consider segments  
$[a_k,u]$, $[v,a_{k+1}]\subset [a_k,a_{k+1}]$ such that $\nu[a_k,u]=\nu[v,a_{k+1}]=\nu_1[C_k]$.
For any $z\in [a_k,v]$ there exists  $z^\circ\le a_{k+1}$ such that 
$\nu[z,z^\circ]=\nu_1[C_k]$. It is easy to see that  
$$
(t\cdot \nu)[a_k,u]\le\nu_1[C_k], \qquad (t\nu)[v,a_{k+1}]\ge (t\cdot\nu_1)[C_k].
$$
Form continuity reasoning there exists
 $[z,z^\circ]$ satisfying the desired property. 
\hfill$\square$

\sm
For each $k$
consider arbitrary measure preserving maps
$$(B_k,\nu)\to (C_k,\nu_1),\qquad (C_k\setminus B_k,\nu_2)\to (C_k,\nu_2).$$
 This produces a map
 $\phi_n$ (see. Fig.\ref{fig:2}). To obtain  
$\psi_n$ we take arbitrary measure preserving maps
$$(B_k,t\cdot\nu)\to (C_k,t\cdot\nu_1),\qquad 
(C_k\setminus B_k,t\cdot \nu_2)\to (C_k,t\cdot\nu_2).$$

Consider a map
$$
\theta_n:\psi_n\circ {\rm id}\circ \phi_n^{-1}:\cL'\to \cL_*'
.$$
The space 
 $\cL'$ consists of two copies  $\R^\times_1$, $\R^\times_2$
of the half-line $\R^\times$, each copy is cutted into segments 
 $C_k$ . The map   $\theta_n$ send each copy of  a segment
 $C_k\subset\R^\times_1$, $C_k\subset\R^\times_2$ to itself, 
 moreover the Radon--Nikodym derivative of
 $\theta_n$ takes values 
$C_k$ in limits $[a_{k},a_{k+1}]$.

It is easy to see that the sequence
 $\theta_n$ converges to the map 
${\rm id}:\cL'\to \cL'_* $.
\hfill $\square$

\sm

{\bf\punct The spreading of measures.}
Denote
$$
\cL:=\cL[\nu,0,\dots,0], \qquad \cL''=\cL[0,0,\dots; \nu].
$$
Let $\cL_*$, $\cL_*''$ be the same measure spaces with the measure multiplied by  $t$.
We construct a sequence of measure preserving bijections 
$$\xi_n: \cL\to\cL'', 
\qquad
\zeta_n: \cL_*\to\cL_*''
.
$$
For this aim, consider the same partitions
 $\frh_n$ of the space $(\R^\times,\nu)$.
 Consider arbitrary measure preserving maps%
\footnote{Recall that any two Lebesgue spaces with  continuous probabilistic measures
are equivalent, see e.g.,  \cite{Roh1}.}
$$
(C_k,\nu)\to (C_k\times [0,1], \nu\times dx),\qquad 
(C_k,t\nu)\to (C_k\times [0,1],(t\nu)\times dx)
$$
This gives us the maps 
 $\xi_n$ и $\zeta_n$. Consider the map 
$$
\upsilon_n=\zeta_n\circ{\rm id}\circ \xi_n^{-1}: \cL''
\to \cL''_*.
$$
The map 
 $\upsilon_n$ sends each   $C_k\times[0,1]$ to itself,
 its Radon--Nikodym derivative on
 $C_k\times[0,1]$ varies in the limits  
$[a_{k-1},a_k]$. Passing to a limit as   $n\to\infty$, 
we get the identity map
 $\cL''\to \cL''_*$.

\sm

{\bf\punct Proof of Theorem   \ref{l:4}}.
Пусть $\nu\in\cM$. Without loss of generality, we can assume that   $\nu$
is continuous. Expand $\nu=\nu_1+\nu_2+\dots+\nu_\infty$.
Set
$$
\cL^k=\cL\bigl[\nu_1,\dots,\nu_k,\sum_{j=k+1}^\infty\nu_j,0,0,\dots;\nu_\infty\bigr],
\qquad \cL^\infty:=\cL[\nu_1,\nu_2,\dots;\nu_\infty].
$$
Let $\cL^k_*$, $\cL_*^\infty$ be the same measure spaces with measures multiplied by   $t$.
Let $\mathrm{id^k}:\cL^k\to \cL^k_*$,  $\mathrm{id^\infty}:\cL^\infty\to \cL^\infty_*$ denote the identical maps.

Iterating arguments of the  two previous subsections, we obtain that the closure of
$S[\nu,0,\dots,0]$ contains  elements  $\mathrm{id}^k$ for any finite  $k$.
Consider a map  $\alpha_k:\cL^k\to\cL^\infty$ constructed in the following way.
It is identical on  $\R_1^\times$, \dots, $\R_k^\times$
send the semi-line $\R^\times_{k+1}$ to $\coprod_{j\ge k+1} \R_j$ preserving the measure.
In the same way we construct a map  $\beta_k:\cL_*^k\to\cL_*^\infty$.
It is easy to see that the sequence
$$
\chi_k:= \beta_k\circ \mathrm{id}_k\circ \alpha_k^{-1}:\cL^\infty\to \cL^\infty_*
.
$$
converges to $\mathrm{id}^\infty$.

\sm

{\bf\punct Construction of the function  $\psi$.} Here we obtain the continuous section
 $\cM\to\Gms$.
 Consider the distribution function
 $z=F(y)$ of the measure $\nu$
 and the inverse function
  $y=G(z)$. If $y_0$ is a discontinuity point of $F$, we set
$G(z)=y_0$ on the segment $[F(y_0-0), F(y_0)+y_0)$. If  $F$  takes some value $z_0$ on
a segment of nonzero length, then   $G(z_0)$ is not defined. Further, we set  
 $\psi(x)=\int_0^x G(z)\,dz$.

\sm

{\bf\punct Proof of Proposition \ref{pr:section}.}
Let $\nu_j$  converges to $\nu$ in $\cM$,
 $y=\psi_j(x)$, $y=\psi(x)$ be the corresponding maps   $[0,1]\to[0,1]$.
We must prove that  $\psi_j$ converges to  $\psi$ in $\Gms$.

\sm

\begin{figure}
$$
\epsfbox{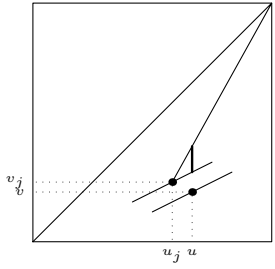}
$$
\caption{To proof of Proposition \ref{pr:section}:
 $(u,v)=H(s)$, $(u_j,v_j)=H_j(s)$.  We mark an interval of possible values of
 $\psi_j(u)$.\label{fig:3}}
\end{figure}

 1) Let $\nu\in\cM$. Consider the map  
$\R^\times\to[0,1]\times[0,1]$ given by the formula 
$$H:s\mapsto\bigl(\nu[(0,s)],(t\cdot\nu)[(0,s)]\bigr).$$
It easy to see that we get the graph of the functions
 $\psi$, from which we remove all straight segments.
The convergence  $\nu_j\to\nu$ means the point-wise convergence of the maps 
 $H_j(s)\to H(s)$. From this it is easy to derive that $\psi_j$ converges to $\psi$
point-wise (See Fig. \ref{fig:3}). In virtue of monotonicity
and continuity of our functions, the point-wise convergence implies the uniform convergence.

\sm

2) Let us show that derivatives  $\psi_j'$ converge  $\psi'$ a.s. Take a point $a$,
where all derivatives 
 $\psi_j'(a)$, $\psi'(a)$ are defined.  Let $\ell_j$,
$\ell$ - be tangent lines to graphs of  $\psi_j$, $\psi$ at $a$.
Suppose that
 $\psi_j'(a)$ does not converge to  $\psi'(a)$. Choose a subsequence 
 $\psi_{n_k}'(a)$  convergent to $\alpha\ne\psi'(a)$. 
 Consider the limit line
 $\ell_{n_k}$, i.e., 
 $$
 \ell^\circ:\qquad y= \alpha(x-a)+\psi(a)
 $$
 It is easy to see
  (for more details, see  \cite{Ale}, Addendum, \S6) that the graph 
  $y=\psi(x)$ is located upper this line. I.e., $\ell^\circ$ is the second
  supporting line at 
 $a$ (the first one was the tangent line),
 this contradicts to the existence of 
  $\psi'(a)$.
 
 \sm
 
3) Now we prove a weak convergence of operators  $T_{1/2}(\psi_j)$ в $L^2[0,1]$.
Let $f$, $h$ be continuous functions.
We must check that the following expressions approach zero
\begin{align}
\Bigl|\int_0^1 f(\psi_j(x))\psi_j'(x)^{1/2}h(x)\,dx -
\int_0^1 f(\psi(x))\psi'(x)^{1/2}h(x)\,dx\Bigr|\le
\nonumber
\\
\le 
\int_0^1 \Bigl|f(\psi_j(x))-f(\psi(x))\Bigr|\,\psi_j'(x)^{1/2}h(x)\,dx
\label{eq:al1}
+\\
+\int_0^1 \Bigl|f(\psi(x))(\psi_j'(x)^{1/2}-\psi'(x)^{1/2})\,h(x)\Bigr|\,dx
\label{eq:al2}
\end{align}

In (\ref{eq:al1}) the convergence   $f(\psi_j(x))\to f(\psi(x))$ is uniform and
$$\int_0^1 \psi_j^{1/2}(x)\le \int_0^1 (\psi_j^{1/2})^2(x) = 1.$$
By the Fatou Lemma, (\ref{eq:al1}) tends to zero. 
Further notice that for functions 
 $\psi\in\cG$ we have a priory estimation 
$$
\psi'(x)\le \frac{1-\psi(x)}{1-x}\le \frac 1{1-x}
.$$
Hence the convergence in the integral
 (\ref{eq:al2}) is dominated on each segment 
$[0,1-\epsilon]$. This implies that integrals  $\int_0^{1-\epsilon}(\dots)$ approach zero.
Further, denote $C=\bigl(\max |f(x)|\cdot \max |g(x)|\bigr)$, 
\begin{multline*}
\int_{1-\epsilon}^1 (\dots)\le C
\int_{1-\epsilon}^1 (\psi_j'(x)^{1/2}+ \psi'(x)^{1/2})\,dx
\le \epsilon C
\int_{1-\epsilon}^1( \psi_j'(x)+ \psi'(x))\,dx
=\\=
\epsilon C\bigl[(1-\psi_j(1-\epsilon)\bigr)+ \bigl(1-\psi(1-\epsilon)\bigr)\bigr]
\end{multline*}
and this value is small for small
 $\epsilon$.
\hfill $\square$


\sm

{\bf\punct Proof of Theorem \ref{th:last}.}
Cut $M$ into pieces $A_{ij}:=M_i\cap g^{-1} M_j$, and also into pieces
$B_{ij}=g M_{ij}=g(M_i)\cap M_j$. We get a collection of maps  
$A_{ij}\to B_{ij}$. Now the question is reduced to a canonical form of each map.

\tt

 Math.Dept., University of Vienna,
 
 Institute for Theoretical and Experimental Physics
 
MechMath Dept., Moscow State University

 neretin(frog)mccme.ru
 
URL: http://www.mat.univie.ac.at/$\sim$neretin/

\end{document}